\documentclass{article}
\usepackage{amsmath,amssymb}
\newtheorem{thm}{Theorem}

\newtheorem{prop}[thm]{Proposition}

\newtheorem{prob}{Problem}
\newtheorem{open}{Open Question}

\newcommand{\bi}{\begin{itemize}}
\newcommand{\ei}{\end{itemize}}
\newcommand{\be}{\begin{enumerate}}
\newcommand{\ee}{\end{enumerate}}
\newcommand{\bc}{\begin{center}}
\newcommand{\ec}{\end{center}}
\newcommand{\bt}{\begin{tabular}}
\newcommand{\et}{\end{tabular}}
\newcommand{\ba}{\begin{array}}
\newcommand{\ea}{\end{array}}

\newcommand{\Z}{\mathbb Z}

  \newcommand{\bp}{\noindent\textit{Proof:}  }
  \newcommand{\ep}{\hfill$\Box$}

\newcommand{\Diek}{Diek}
\newcommand{\EFR}{EFR}
\newcommand{\Mias}{MR2645045}
\newcommand{\GARM}{MR2525867}
\newcommand{\WordProc}{MR1161694}
\newcommand{\NS}{MR1329042}
\newcommand{\Groves}{MR1425318}
\newcommand{\Fordham}{MR1998934}
\newcommand{\Dehn}{MR1511645}
\newcommand{\Collins}{MR840121}
\newcommand{\BelkBrown}{MR2197808}
\newcommand{\Elder}{Elder}

\title{Some geodesic problems in groups}

\author{
Murray Elder\\
\small Mathematics \\[-0.8ex]
\small University of Queensland, Brisbane, Australia \\[-0.8ex]
\small\texttt{murrayelder@gmail.com}\\[1.2ex]
Andrew Rechnitzer\\
\small Department of Mathematics,  \\[-0.8ex]
\small University of British Columbia, Vancouver, Canada\\[-0.8ex]
\small\texttt{andrewr@math.ubc.ca}
}

\date{\small \today. }

\begin{document}
\maketitle

\begin{abstract}
We consider several algorithmic problems concerning  geodesics in finitely generated groups. We show that the three geodesic problems considered by Miasnikov {\em et al} 
are polynomial-time reducible to each other. We study two new geodesic problems which arise in a previous paper of the authors and Fusy. 
\end{abstract}

\section{Introduction} 
The study of algorithmic problems in group theory goes back to Dehn \cite{\Dehn}, and has been a key theme in much of combinatorial and geometric group theory since. In a recent paper Miasnikov {\em et al} \cite{\Mias}  considered three algorithmic problems concerning geodesics in finitely generated groups. In this article we show that these three problems are in fact equivalent, in the sense that each is polynomial time reducible to the others. We  consider two more related problems which arose in work of the authors with \'Eric Fusy, in computing the growth and geodesic growth rates of groups \cite{\EFR}. We show that these new problems have efficient polynomial time solutions for a large class of groups. 

The two problems concern how adding a generator to a  geodesic effects  length. Let $\langle \mathcal G\;|\;\mathcal R\rangle$ be a presentation for a group $G$, with $\mathcal G$ finite, and  $\ell(u)$  the length of a geodesic
representative for a word $u\in \left(\mathcal G^{\pm 1}\right)^*$. 
\begin{prob}
Given a geodesic $u$ and generator $x\in\mathcal G$, decide whether $\ell(ux)-\ell(u)=0,1$ or $-1$.
\end{prob}
Note that if all the relators in $\mathcal R$ have even length, then   $\ell(ux)=\ell(u)$ would imply a word of odd length equal to the identity, so $\ell(ux)-\ell(u)=0$ is not an option. 
 So the problem turns into the decision problem:
\begin{prob}
Given a geodesic $u$ and generator $x\in\mathcal G$, decide if $\ell(ux)>\ell(u)$.
\end{prob}

In \cite{\EFR} we use the fact that Thompson's group $F$ has an efficient (time $O(n^2\log n)$ and space  $O(n\log n)$) solution to Problem 1, together with a technique from \cite{\GARM}, to compute the number of geodesics and elements up to length 22. We remark that the same procedure will work efficiently for any group and generating set that has an efficient (polynomial space and time) solution to Problem 1.

This article is organised as follows. In Section \ref{sec:Mias} we consider the geodesic problems of Miasnikov {\em et al} \cite{\Mias}. We prove that they are polynomially reducible  to each other, and that Problems 1 and 2 polynomially reduce to them.
In Section \ref{sec:WP} we prove that if $\mathcal R$ is countably enumerable, the problems of Miasnikov {\em et al} reduce to Problems 1 and 2, this time not necessarily in polynomial time and space. This implies a solvable word problem, and so we have examples (from \cite{\Collins} for instance) for which Problems 1 and 2 are unsolvable. In Section \ref{sec:egs} we describe groups which have polynomial time solutions to Problems 1 and 2 and give some open problems.

\textit{Acknowledgments} We thank Volker Diekert and Andrew Duncan for explaining how Problems 1 and 2 imply a solvable word problem and geodesic problem for recursively presented groups. Thanks to Alexei Miasnikov, Vladimir Shpilrain and Sasha Ushakov for many fruitful discussions about this work, and to the anonymous reviewer for their careful proofreading and suggestions.

\section{Geodesic problems of Miasnikov {\em et al}}\label{sec:Mias}

In \cite{\Mias} Miasnikov {\em et al} consider the following algorithmic and decision problems for a group $G$ with finite generating set $\mathcal G$.
\begin{prob}
(Geodesic problem) Given a  word in $\mathcal G^{\pm 1}$, find a geodesic representative for it.
\end{prob}
\begin{prob}
(Geodesic length problem) Given a  word in  $\mathcal G^{\pm 1}$, find the length of a geodesic representative.
\end{prob}
\begin{prob}
(Bounded geodesic length problem) Given a  word in $\mathcal G^{\pm 1}$ and an integer $k$, decide if a geodesic representative has length $\leq k$.
\end{prob}

They show that for free metabelian groups (with standard generating sets), Problem 5 is NP-complete. They also show that
a polynomial time solution to any of these problems implies a polynomial time solution to the next, and each implies a polynomial time solution to the word problem for the group.

\begin{prop}
If $G$ is a group with  finite generating set $\mathcal G$, then Problems 3-5 of Miasnikov {\em et al} are polynomial time and space reducible to each other.\end{prop}
\bp
It is clear that a solution to Problem 3 solves Problem 4,  and a solution to Problem 4 solves Problem 5 . Suppose we can solve Problem 4 in time $f(n)$ and space $g(n)$.
Then consider the following solution to Problem 3. Given a word $u$ of length $n$ in the finite (inverse closed) generating set $\mathcal G$, with $|\mathcal G|=k$, apply the solution to Problem 4 to  $u$. If it returns  $n$, then $u$ is a geodesic. Else the output is $m<n$. Pick a generator $x$ and run Problem 4 on $ux$. If the output is not $m-1$, then pick another $x$. If the output is $m-1$, set $x_1=x$. For $i=2,\ldots, m$, pick a generator $x$, and run Problem 4 on $ux_1\ldots x_{i-1}$. If output is $m-i$ then set $x_i=x$, else pick again. After $m$ iterations we have $ux_1x_2\ldots x_m$ has length 0 so is equal to the identity, so $x_m^{-1}\ldots x_2^{-1}x_1^{-1}$ is a geodesic for $u$.
Each iteration takes at most $kf(n+i)$ time, and $i\leq m<n$, 
so in total this algorithm takes $nkf(2n)$ time, and space $g(n+m)$.

Next,  suppose we can solve the Problem 5 in time $f(n)$ and space $g(n)$.
Then consider the following solution to Problem 4. Given a word $u$ of length $n$, 
run the solution to Problem 5 on the pair $(u,n-1)$.  If the output is No, then $u$ is a geodesic, so output $n$. While the answer is Yes, run Problem 5 on  $(u,n-i)$ for $i=2,3,\ldots$ until the answer is No, and thus the length of a geodesic for $u$ is   $n-i$. The total time is at most $nf(n)$ and space $g(n)$.
\ep

This answers Problem 5.3 in \cite{\Mias}.

\begin{prop}
A polynomial time and space solution to Problem 5 implies a polynomial time and space solution to Problem 1.
\end{prop}
\bp
Given a geodesic word $u=u_1\ldots u_{\ell(u)}$ and a generator $x\in\mathcal G$, run  Problem 5 on  $(ux,\ell(u)-1)$. If it returns Yes, then $\ell(ux)-\ell(u)=-1$. If not, then run Problem 5 on $(ux,\ell(u))$. If this returns Yes, $\ell(ux)-\ell(u)=0$. If not, then $\ell(ux)-\ell(u)=1$.
\ep

The converse of this proposition, that a polynomial time and space solution to Problem  2 implies a polynomial solution to Problems 3-5, is not obvious. In the next section we prove that the solvability of Problem 2 implies some solution to Problems 3-5, but not preserving time and space complexity.

\section{Word problem}\label{sec:WP}

If $G=\langle \mathcal G\;|\; \mathcal R\rangle$ is recursively presented, by which we mean $\mathcal R$  is countably enumerable, and has a solution to Problem 1 or 2, then a very brute force procedure which runs through all possible words can solve the geodesic problem (Problem 3). It follows that Problems 4 and 5, and the word problem, are solvable. Since there are many groups with unsolvable word problem, including finitely presented examples \cite{\Collins}, then there are groups for which  Problems 1 and 2 are unsolvable. 
This also shows that a solution to Problem 2 implies some solution to Problem 1 for any recursively presented group.

The following proof  was described to us by Volker Diekert and Andrew Duncan. The notation $u=_Gv$ means that $u$ and $v$ represent the same group element.

\begin{prop}\label{prop:WP}
If $G=\langle \mathcal G \; | \; \mathcal R\rangle$ with $\mathcal G$ finite and $R$  countably enumerable, and has a solution to Problem 2, then Problem 3 is solvable. 
\end{prop}

\bp
We proceed by induction on the length of the input word $u$ to Problem 3. 
If $|u|=0$ then clearly $u$ is geodesic. Assume for all words of length $n\geq 0$ we can find a geodesic representative, and consider the word $w$ of length  $n+1$. Write $w=ux$ where $x$ is a generator. So we can find a geodesic representative $v$ for $u$ by inductive assumption.

Input $v,x$ into Problem 2. If it returns $\ell(vx)>\ell(x)$, then $vx$ is a geodesic for $w$. 
If not, then we know that $vx=_G z$ for some word $z$ of length $\leq n$. 
So $1=_G z(vx)^{-1}=_Gc_1 c_2 .... c_k$, for some $k$, where each $c_i$ is a conjugate of an element of $\mathcal R$,
and so $z=_G c_1c_2\ldots c_k(vx)$.

Since $\mathcal R$ is countably enumerable, so too is the set of conjugates, as is the set of products of conjugates. If the list of products of conjugates is $p_1,p_2,\ldots$ then running through $p_1(vx), p_2(vx),\ldots$, freely reducing each one, we must eventually find a word of freely reduced length $\leq n$ equal to $z=_Gvx$. 

If this word has length $n-1$ then it must be a geodesic for $vx$ (since $v$ was geodesic of length $n$). Otherwise it has length $n$, so applying the inductive assumption we can find a geodesic for it.
\ep

\section{Examples}\label{sec:egs}

We end by considering groups that have efficient solutions to Problems 1 and 2.
As noted in the introduction, Thompson's group $F$ has an efficient solution to Problem 4 (and thus Problems 1-5). This result is essentially due to Fordham \cite{\Fordham}, who first introduced a simple technique for computing the geodesic length of an element represented as a pair of rooted binary trees. The implementation used in \cite{\EFR} which gives the $O(n^2\log n)$ time and  $O(n\log n)$ space solution is the version of Fordham's technique due to Belk and Brown \cite{\BelkBrown}. So what other groups have polynomial time and space solutions to Problems 1 and 2?

\begin{prop}
If the full set of geodesics for a group with a finite generating set forms a regular language, and one can construct the corresponding finite state automaton, then Problem 2 is solvable in linear time and constant space.
\end{prop}
\bp
Given a word $u$ and a generator $x$, read $ux$ into  the finite state automaton.  $ux$ is accepted if and only if $\ell(ux)>\ell(u)$. The space required is proportional to the number of states in  the automaton.
\ep

It follows that all groups satisfying the  {\em falsification by fellow traveler property} \cite{\NS} have a linear time solution to Problem 2, provided the fellow traveling constant is known. Examples include hyperbolic groups, abelian groups, Coxeter groups, and virtually abelian groups with some generating sets  \cite{\NS}. 

An automatic structure on a group is {\em strongly geodesic}  if it includes the set of all geodesics \cite{\WordProc}. In this case our argument to solve Problem 2 may not work -- if $ux$ is rejected by the acceptor automaton then we know $\ell(ux)\leq\ell(u)$, but non-geodesic words may also be accepted. More generally, we may ask:
\begin{open}
Does every automatic group have a polynomial time and space solution to  Problems 1 and 2?
\end{open}


In  \cite{\Elder} the first author gives an algorithm to compute a geodesic representative for a word in the solvable Baumslag-Solitar groups $\langle a,t \; | \; tat^{-1}=a^n\rangle$ in linear time and space, thus solving Problem 3. By \cite{\Groves} the full set of geodesics for these groups (with this generating set) fails to be regular.
An intriguing open problem is the following:
\begin{open}
For the non-solvable Baumslag-Solitar groups $\langle a,t \; | \; ta^mt^{-1}=a^n\rangle$ with $|m|,|n|>1$, is there an algorithm to solve any of Problems 1-5 in polynomial time and space?
\end{open}
In the case $m=\pm n$ the group is isomorphic to $F_{|n|}\times \Z$ 
so the answer is yes (such groups enjoy the falsification by fellow traveler property, for example). Recently  Volker Diekert and J\"urn Laun
partially answered this question -- they give a clever geodesic normal form which yields a polynomial solution for the cases where $n$ divides $m$ \cite{\Diek}.

As noted, for presentations which have only even length relators, Problems 1 and 2 are solved by the same algorithm. For a presentation with odd length relators, it is not clear that  a solution to Problem 2 solves Problem 1 in the same time and space complexity.

\begin{open} Is there an example with odd length relators where Problems 1 and 2 do not have the same time and space complexities?
\end{open}

And lastly
\begin{open} Is there an example where Problems 1 and/or 2 have polynomial time (and space) solutions, but problems 3-5 are superpolynomial?
\end{open}

\bibliographystyle{plain}
\bibliography{geodprobs-refs} 

\end{document}